\documentclass[12pt]{amsart}
\usepackage{latexsym}
\usepackage{amsthm}
\usepackage{amsmath}
\usepackage{amsfonts}
\usepackage{amssymb}
\usepackage[dvips]{graphicx}
\usepackage{xypic}
\input xy
\xyoption{all}

\newtheorem{theo}{Theorem}[section]
\newtheorem{lemma}[theo]{Lemma}
\newtheorem{propo}[theo]{Proposition}
\newtheorem{defi}[theo]{Definition}
\newtheorem{coro}[theo]{Corollary}
\newtheorem{rem}[theo]{Remark}
\newtheorem{exam}[theo]{Example}

\newcommand\Ab{\operatorname{\bf Ab}}

\newcommand\Iso{\operatorname{Iso}}

\newcommand\card{\operatorname{card}}
\newcommand\id{\operatorname{id}}

\newcommand\Top{\operatorname{\bf Top}}
\newcommand\Pos{\operatorname{\bf Pos}}

\newcommand\cof{\operatorname{cof}}

\newcommand\ca{\mathcal {A}}

\newcommand\cc{\mathcal {C}}
\newcommand\cd{\mathcal {D}}
\newcommand\cf{\mathcal {F}}
\newcommand\ch{\mathcal {H}}

\newcommand\ce{\mathcal {E}}

\newcommand\ck{\mathcal {K}}
\newcommand\cl{\mathcal {L}}
\newcommand\cm{\mathcal {M}}
\newcommand\crr{\mathcal {R}}

\newcommand\cw{\mathcal {W}}
\newcommand\cx{\mathcal {X}}

\date{August 11, 2007}
 
\begin{document}
\title{On combinatorial model categories}
\author[J. Rosick\'{y}]
{J. Rosick\'{y}$^*$}
\thanks{ $^*$ Supported by MSM 0021622409 and GA\v CR 201/06/0664.} 
\address{\newline J. Rosick\'{y}\newline
Department of Mathematics and Statistics\newline
Masaryk University, Faculty of Sciences\newline
Jan\'{a}\v{c}kovo n\'{a}m. 2a, 60200 Brno, Czech Republic\newline
rosicky@math.muni.cz
}
\begin{abstract}
Combinatorial model categories were introduced by J. H. Smith as model categories which are locally presentable 
and cofibrantly generated. He has not published his results yet but proofs of some of them were presented by T. Beke 
or D. Dugger. We are contributing to this endeavour by proving that weak equivalences in a combinatorial model
category form an accessible category. We also present some new results about weak equivalences and cofibrations 
in combinatorial model categories.
\end{abstract}

\maketitle
 
\section{Introduction}
Model categories were introduced by Quillen \cite{Q} as a foundation of homotopy theory. Their modern theory
can found in \cite{Ho} or \cite{H}. Combinatorial model categories were introduced by J. H. Smith as model categories 
which are locally presentable and cofibrantly generated. The latter means that both cofibrations and trivial cofibrations 
are cofibrantly generated by a set of morphisms. He has not published his results yet but some of them can be found 
in \cite{B} or \cite{D}. In particular, \cite{B} contains the proof of the theorem characterizing when a class $\cw$ 
of weak equivalences makes a locally presentable category $\ck$ to be a combinatorial model category with a given 
cofibrantly generated class $\cc$ of cofibrations. The characterization combines closure properties of $\cw$ together 
with a smallness condition saying that $\cw$ satisfies the solution set condition at the generating set $\cx$ 
of cofibrations. We will show that these conditions are also necessary. This is based or another result of J. H. Smith 
saying that, in a combinatorial model category, $\cw$ is always accessible and accessibly embedded in the category $\ck^\to$
of morphisms of $\ck$ (he informed me about this result in 2002 without indicating a proof). Since \cite{D} proved
that $\cw$ is accessibly embedded into $\ck^\to$, it remains to show that $\cw$ is accessible. We will provide
the mi\-ssing proof. Our argument is based on the fact that homotopy equivalences form a full image of an accessible
functor into $\ck^\to$. Surprisingly, the smallness condition above is automatic in a set theory with a strong axiom 
of infinity called Vop\v enka's principle. We will show that Vop\v enka's principle is equivalent to the fact that 
the smallness condition can be avoided. 

Nearly all model categories are either combinatorial or Quillen equivalent to a combinatorial one. The example
of the latter are topological spaces because $\Top$ is not locally presentable (although both cofibrations and
trivial cofibrations are cofibrantly generated). $\Top$ is Quillen equivalent to the combinatorial model category
of simplicial sets. Let us mention that J. H. Smith claimed that one can make $\Top$ a combinatorial model
category by reducing topological spaces to simplex-generated ones (we were able to prove this statement in \cite{FR}).
But there are also model categories which are neither locally presentable nor cofibrantly generated (see \cite{CD})
or model categories which are locally presentable but not cofibrantly generated (see \cite{AHRT1} or \cite{CH}).
All known combinatorial model categories have the class $\cc$ of cofibrations accessible and accessibly embedded
into $\ck^\to$. The consequence is that the full subcategory $\ck_{cf}$ of $\ck$ consisting of cofibrant and fibrant 
objects is accessible and accessibly embedded into $\ck$. We will show that it does not need to be true in general,
however, one can prove that both $\cc$ and $\ck_{cf}$ are the closure under retracts of a full image of an accessible
functor (into $\ck^\to$ or into $\ck$). This is used in the above mentioned proof of accessibility of $\cw$ because
a morphism $f$ is a weak equivalence if and only if its replacement $R(f)$ to $\ck_{cf}$ is a homotopy equivalence.

There remains one claim of J. H. Smith which I have not not able to prove. It concerns the existence of the smallest
$\cw$ making a locally presentable category $\ck$ a combinatorial model category with a given cofibrantly generated 
class $\cc$ of cofibrations. This is true under Vop\v enka's principle but I am not able to avoid its use.  

\section{Accessible categories and their generalizations}
The theory of accessible categories were created by M. Makkai and R. Par\' e \cite{MP}. They departed from
\cite{GU} and \cite{AGV} and their motivation was model theoretic (these models come from logic and not from
homotopy theory). Applications of this theory to homotopy theory have appeared during the last ten years. 
I add that, a month ago, G. Maltsiniotis sent me the preliminary redaction of the unpublished Grothendieck's 
manuscript \cite{G} where he independently develops a theory of accessible categories on his own, motivated 
by homotopy theory in this case. In \cite{AR}, we showed that quite a few of properties of accessible categories 
depend on set theory, in particular on Vop\v enka's principle. This is a large cardinal axiom implying the existence 
of a proper class of measurable (compact, or even extendible) cardinals. On the other hand, its consistency follows 
from the exiatence of a huge cardinal. One of benefits of using accessible categories in homotopy theory is 
the realization that some open homotopy theoretic problem depend (or may depend) on set theory (see \cite{CSS}, 
\cite{C1}, \cite{M} or \cite{CGR}). 
 
Let us recall that a category $\ck$ is called $\lambda$-\textit{accessible}, where $\lambda$ is
a regular cardinal, provided that
\begin{enumerate}
\item[(1)] $\ck$ has $\lambda$-directed colimits,
\item[(2)] $\ck$ has a set $\ca$ of $\lambda$-presentable objects such that every object of $\ck$ is a 
$\lambda$-directed colimit of objects from $\ca$.
\end{enumerate}

A category is \textit{accessible} if it is $\lambda$-accessible for some regular cardinal $\lambda$. A cocomplete 
accessible category is called \textit{locally presentable}. All needed facts about locally presentable and accessible 
categories can be found in \cite{AR} or \cite{MP}. A full subcategory $\cl$ of an accessible category $\ck$ is called
\textit{accessibly embedded} if it is closed under $\lambda$-directed colimits for some regular cardinal $\lambda$.

\begin{propo}\label{prop2.1}
Let $\ck$ be an accessible category. Any union of a set of accessible and accessibly embedded subcategories of $\ck$
is accessible and accessibly embedded in $\ck$.
\end{propo}
\begin{proof}
By \cite{AR}, 2.36, a full subcategory of $\ck$ is accessible and accessibly embedded if and only if it is closed
under $\lambda$-directed colimits and $\lambda$-pure subobjects for some regular cardinal $\lambda$. Thus, given a set
$\cl_i$, $i\in I$ of accessible and accessibly embedded subcategories of $\ck$, there is a regular cardinal $\lambda$
such that all $\cl_i$, $i\in I$ are closed under $\lambda$-directed colimits and $\lambda$-pure subobjects. Then their
union $\cl$ is closed under $\lambda$-pure subobjects as well. Without loss of generality, we can assume that
$\card I<\lambda$. Consider a $\lambda$-directed diagram $D:\cd\to\cl$. Let $D_i:\cd_i\to\cl_i$ be the pullback
of $D$ with the inclusion $\cl_i\to\cl$. There is $i\in I$ such that the embedding of $\cd_i$ to $\cd$ is cofinal. 
In fact, assuming the contrary, there is an object $d_i$ in $\cd$ with $d_i\downarrow\cd_i=\emptyset$, for each
$i\in I$, which contradicts to $\card I<\lambda$. Consequently, $\cl$ is closed under $\lambda$-directed colimits 
and thus it is accessible and accessibly embedded.
\end{proof}

Let $F:\cl\to\ck$ be an accessible functor. Recall that it means that both $\cl$ and $\ck$ are accessible and $F$
preserves $\lambda$-directed colimits for some regular cardinal $\lambda$. The full subcategory of $\ck$ consisting 
of objects $FL$, $L\in\cl$ is called a \textit{full image} of $F$. 
While accessible categories are, up to equivalence,
precisely categories of models of basic theories, full images of accessible functors are, up to equivalence, precisely
categories of structures which can be axiomatized using additional operation and relation symbols (see \cite{R1});
they are also called pseudoaxiomatizable. In both cases, we use infinitary first-order theories. In \cite{AR1}, 
a category $\ck$ was called \textit{preaccessible} if there is a regular cardinal $\lambda$ such that $\ck$ has 
a set $\ca$ of $\lambda$-presentable objects such that every object in $\ck$ is a $\lambda$-directed colimit 
of objects from $\ca$. Hence accessible categories are precisely preaccessible ones having $\mu$-directed colimits 
for some regular cardinal $\mu$.

\begin{propo}\label{prop2.2}
Assuming the existence of a proper class of compact cardinals, every full image of an accessible functor is preaccessible.
\end{propo}
\begin{proof}
Following the uniformization theorem (see \cite{AR}, 2.19), there is a regular cardinal $\mu$ such that both $\cl$ 
and $\ck$ are $\mu$-accessible and $F$ preserves $\mu$-directed colimits and $\mu$-presentable objects. There is 
a compact cardinal $\lambda>\mu$ and, by \cite{R1}, Theorem 1, the inclusion of $F(\cl)$ to $\ck$ preserves 
$\mu$-directed colimits. Since a compact cardinal is (strongly) inaccessible, we have $\lambda\trianglelefteq\mu$ 
(see \cite{AR}, 2.12) and thus $F$ preserves $\mu$-presentable objects (see \cite{AR}, 2.18 (2)). Hence $FL$ is 
$\mu$-presentable in $F(\cl)$ for each $\mu$-presentable object $L$ in $\cl$.
\end{proof}

\begin{rem}\label{rem2.3}
{
\em
Let $\cl\to\ck$ be an accessible functor and assume that its full image $F(\cl)$ is closed under $\lambda$-directed
colimits in $\ck$ for some regular cardinal $\lambda$. Then, in the same way as in the proof of \ref{prop2.2},
we show that $F(\cl)$ is accessible.
}
\end{rem}

Let $\cl$ be a full subcategory of a category $\ck$ and $K$ an object in $\ck$. We say that $\cl$ satisfies 
the \textit{solution-set condition} at $K$ if there exists a set of morphisms $(K\to L_i)_{i\in I}$ with $L_i$ 
in $\cl$ for each $i\in I$ such that every morphism $f:K\to L$ with $L$ in $\cl$ factorizes through some $f_i$, 
i.e., $f=gf_i$. $\cl$ is called \textit{cone-reflective} in $\ck$ if it satisfies the solution-set condition
at each object $K$ in $\ck$ (see \cite{AR}). Given a set $\cx$ of objects of $\ck$, we say that $\cl$ 
satisfies the solution set condition at $\cx$ if it satisfies this condition at each $X\in\cx$.  

\begin{propo}\label{prop2.4}
The full image of an accessible functor $F:\cl\to\ck$ is cone-reflective in $\ck$.
\end{propo}
\begin{proof}
Let $K$ be an object in $\ck$. Like in \ref{prop2.2}, we can assume that there is a regular cardinal $\mu$ such that
$K$ is $\mu$-presentable, both $\cl$ and $\ck$ are $\mu$-accessible and  $F$ preserves $\mu$-directed colimits and 
$\mu$-presentable objects. Hence each morphism $f:K\to F(L)$ factorizes through some $K\to F(L')$ with $L'$
$\mu$-presentable in $\cl$. It yields the solution set condition at $K$.
\end{proof}

Recall that an idempotent $f:K\to K$ in a category $\ck$ (i.e., $ff=f$) \textit{splits} if there exist morphisms
$i:X\to K$ and $p:K\to X$ such that $pi=\id_X$ and $ip=f$. In an accessible category, all idempotents split. Every 
category $\ck$ has a \textit{split idempotent completion} $\tilde{\ck}$. This means that there is a functor 
$H:\ck\to\tilde{\ck}$ where all idempotents split in $\tilde{\ck}$ such that, whenever $G:\ck\to\cm$ is a functor 
into $\cm$ where all idempotents split, than there is a unique (up to an isomorphism) functor $G':\tilde{\ck}\to\cm$ 
with $G'H\cong G$ (see \cite{AR}, Ex. 2.b).

\begin{lemma}\label{lem2.5}
Let $\cl$ be a cone-reflective subcategory of an accessible category $\ck$. Then the split idempotent completion
$\tilde{\cl}$ is cone-reflective in $\ck$.
\end{lemma}
\begin{proof}
Since all idempotents split in $\ck$, $\cl$ is a full subcategory of $\ck$. In fact, $\tilde{\cl}$ is the closure
of $\cl$ under retracts in $\ck$. Let $f_i:K\to L_i$, $i\in I$ be a cone-reflection of $K$ to $\cl$. 
For each subobject $i:X\to L$, $L\in\cl$, split by $p:L\to X$ and each $f:K\to X$, there is $g:L_i\to L$ with
$gf_i=if$. Hence $f=pif=pgf_i$. Thus $f_i$ is a cone-reflection of $K$ to $\tilde{\cl}$.
\end{proof} 

I do not know whether a split idempotent completion of a full image of an accessible functor is a full image 
of an accessible functor.

The powerfull theorem of Makkai and Par\' e says that accessible categories are closed under all pseudolimits
(see \cite{MP}, 5.1.6 or \cite{AR}, Exercise 2.n). In particular, they are closed under pseudopullbacks. Recall
that pseudo means commutativity up to isomorphism. It is easy to see that this result can be extended to full images 
of accessible functors.

\begin{lemma}\label{lem2.6}
Let $F_1:\cm_1\to\ck$ and $F_2:\cm_2\to\ck$ be accessible functors and $\cl_1,\cl_2$ their full images. Let
$\cl$ be a pseudopullback
$$
\xymatrix@=4pc{
\cl_1 \ar[r]^{} & \ck \\
\cl \ar [u]^{} \ar [r]_{} &
\cl_2 \ar[u]_{}
}
$$
Then $\cl$ is a full image of an accessible functor $\cm\to\ck$.
\end{lemma}
\begin{proof}
It suffices to take $\cm$ as a pseudopullback
$$
\xymatrix@=4pc{
\cm_1 \ar[r]^{R} & \ck \\
\cl \ar [u]^{} \ar [r]_{} &
\cm_2 \ar[u]_{}
}
$$
Then $\cl$ is the full image of the induced functor $F:\cm\to\ck$.
\end{proof}

\section{Weak factorization systems}
 
Let $\ck$ be a category and $f: A\to B$, $g: C\to D$ morphisms
such that in each commutative square
$$
\xymatrix@=4pc{
A \ar[r]^{u} \ar[d]_{f}& C \ar[d]^g\\
B\ar[r]_v & D
}
$$
there is a diagonal $d:B \to C$ with $df=u$ and $gd=v$.
Then we say that $g$ has the \textit{right lifting property}
w.r.t. $f$ and $f$ has the \textit{left lifting property} w.r.t.
$g$. 
For a class $\cx$ of morphisms of $\ck$ we put
\begin{align*}
\cx^{\square}& = \{g|g \ \mbox{has the right lifting property
w.r.t.\ each $f\in \cx$\} and}\\
{}^\square\cx & = \{ f|f \ \mbox{has the left lifting property
w.r.t.\ each $g\in \cx$\}.}
\end{align*}

\begin{defi}\label{def3.1}
{\em    
A \textit{weak factorization system} $(\cl,\crr)$ in a
category $\ck$ consists of two classes $\cl$ and $\crr$ of morphisms
of $\ck$ such that
\begin{enumerate}
\item[(1)] $\crr = \cl^{\square}$, $\cl = {}^\square \crr$, and
\item[(2)] any morphism $h$ of $\ck$ has a factorization $h=gf$ with
$f\in \cl$ and $g\in \crr$.
\end{enumerate}

A weak factorization system $(\cl,\crr)$ is called \textit{cofibrantly generated} if there is a set
$\cx$ of morphisms such that $\crr=\cx^\square$.
}
\end{defi}

This definition and the following basic facts can be found in \cite{B} (or \cite{AHRT}).

\begin{rem}\label{rem3.2}
{
\em
(1) Given a weak factorization system $(\cl,\crr)$ then $\cl$ is 
\textit{cofibrantly closed} in the sense
that it contains all isomorphisms and is

(a) stable under pushout,

(b) closed under transfinite composition, and

(c) closed under retracts in comma categories $\ck^\to$.

\noindent
The first condition says that if 
$$
\xymatrix@=4pc{
B \ar[r]^{\overline g}& D\\
A \ar[u]^{f} \ar[r]_{g}& C \ar[u]_{\overline{f}}
}
$$ 
is a pushout and $f\in\cl$ then $\overline f\in\cl$. The second condition means that $\cl$ is closed under composition 
and  if $(f_{ij}:A_i \to A_j)_{i\leq j\leq \lambda}$ is a \textit{smooth} chain (i.e., $\lambda$ is a limit ordinal, 
$(f_{ij}:A_i \to A_j)_{i<j}$ is a colimit for any limit ordinal $j\leq\lambda$) and $f_{ij}\in\cl$ for each 
$i\leq j< \lambda$ then $f_{0\lambda}\in\cl$. In the third condition, $\ck^\to$ denotes the category of morphisms
of $\ck$.

In what follows, $\cof(\cx)$, denotes the cofibrant closure of $\cx$, i.e., the closure of $\cx$ under isomorphisms
and constructions (a)--(c). We always have
$$
\cof(\cx)\subseteq {}^\square(\cx^\square).
$$

(2) Let $\ck$ be a locally presentable category. Then each set $\cx$ of morphisms determines a (cofibrantly generated)
weak factorization system $(\cof(\cx),\cx^\square)$. 
}
\end{rem}

\begin{propo}\label{prop3.3}
Let $\ck$ be a locally presentable category and $\cx$ a set of morphisms. Then $\cx^\square$ is an accessible
category which is accessibly embedded in $\ck^\to$.
\end{propo}

\begin{proof}
It suffices to observe that $g$ has a right lifting property w.r.t. 
$f:A\to B$ if and only if $g$ is injective 
in $\ck^\to$ to the morphism
$$
(f,\id_B):f\to\id_B.
$$
The result then follows from accessibility of small-injectivity classes (see \cite{AR}, 4.7).
\end{proof}

\begin{propo}\label{prop3.4}
Let $\ck$ be a locally presentable category and $\cx$ a set of morphisms.
 Then $\cof(\cx)$ is a split idempotent
completion of a full image of an accessible functor $\cm\to\ck^\to$.
\end{propo}

\begin{proof}
By the proof of \cite{R}, 3.1, there is an accessible functor 
$$
F:\ck^\to\to\ck^\to
$$ 
such that a $(\cof(\cx),\cx^\square)$ factorization of a morphism $f:A\to B$ is $gF(f)$. Since $f$ belongs 
to $\cof(\cx)$ if and only if it is a retract of $F(f)$, $\cof(\cx)$ is a split idempotent completion of the image 
of $F$.
\end{proof}

The following examples show that $\cof(\cx)$ is not always accessible. In these examples, it is a full image
of an accessible functor. I do not know any example where it is not so.

\begin{exam}\label{ex.3.5}
{
\em
(1) Let $\Pos$ be the category of posets and let $\cl$ consist of split monomorphisms. It is easy to see
that $\cl$ is the cofibrant closure of split monomorphisms between finite posets. Since $\Pos$ is locally finitely
presentable, the closure of split monomorphisms under $\lambda$-directed colimits in $\Pos^\to$ precisely consists 
of $\lambda$-pure monomorphisms (see \cite{AR}, 2.40). It is easy to see that, for each regular cardinal $\lambda$, 
there is a $\lambda$-pure monomorphism which does not split. On the other hand, $\cl$ is the full image 
of an accessible functor $F:\cm\to\ck^\to$ where objects of $\cm$ are pairs $(i,p)$ with $pi=\id$ and morphisms
are pairs of morphisms $(u,v):(i,p)\to (i',p')$ such that both $vi=i'u$ and $us=s'v$.

(2) Let $\Ab$ be the category of abelian groups and let $\cl$ be cofibrantly generated by the morphism
$$
0\to\bf Z.
$$
Then the comma category $0\downarrow\cl$ precisely consists of morphisms $0\to A$ such that $A$ is free.
If $\cl$ is accessible then $0\downarrow\cl$ is accessible as well (see \cite{AR}, 2.44). Hence the category $\cf$
of free abelian groups is accessible. However, this statement depends on set theory: (i) $\cf$ is accessible
if there is a compact cardinal and (ii) if the axiom of constructibility is assumed, then $\cf$ is not accessible
(see \cite{MP}, 5.5). 
}

\end{exam}

\begin{coro}\label{cor3.6} 
Let $\ck$ be a locally presentable category and $\cx$ a set of morphisms. Then $\cx^\square$ is cone-reflective
in $\ck^\to$.
\end{coro}

\begin{proof}
It follows from \ref{prop3.4}, \ref{prop2.4} and \ref{lem2.5}.
\end{proof}

\begin{defi}\label{def3.7}
{
\em
Let $\ck$ be a category with finite coproducts equipped with a weak factorization system $(\cl,\crr)$. 
A \textit{cylinder object} $C(K)$ of an object $K$ is given by an $(\cl,\crr)$ factorization of the codiagonal
$$
\nabla : K+K \xrightarrow{\quad  \gamma_K\quad} C(K)
             \xrightarrow{\quad \sigma_K\quad} K
$$
We denote by
$$
\gamma_{1K},\gamma_{2K}:K\to C(K)
$$ 
the compositions of $\gamma_K$ with the coproduct injections. 
}
\end{defi}

This definition was suggested in \cite{KR}. As usual, we say that morphisms $f,g:K\to L$ are \textit{homotopic},
and write $f\sim g$, if there is a morphism $h:C(K)\to L$ such that the following diagram commutes
$$
\xymatrix@=3pc{
K+K \ar[rr]^{(f,g)}
\ar[dr]_{\gamma_K} && L\\
& C(K) \ar[ur]_h
}
$$
Here, $(f,g)$ is induced by $f$ and $g$. The homotopy relation  $\sim$ is clearly reflexive, symmetric, 
compatible with the composition and does not depend on the choice of a cylinder object. But, it is not transitive 
in general and we will denote its transitive hull by $\approx$. We get the quotient functor
$$
Q:\ck\to\ck/\approx.
$$

A morphism $f:K\to L$ is called a \textit{homotopy equivalence} if $Qf$ is the isomorphism, i.e., if there exists
$g:L\to K$ such that both $fg\approx\id_L$ and $gf\approx\id_K$. The full subcategory of $\ck^\to$ consisting
of homotopy equivalences w.r.t. a weak factorization system $(\cl,\crr)$ will be denoted by $\ch_\cl$. 

\begin{propo}\label{prop3.8}
Let $\ck$ be a locally presentable category and $\cx$ a set of morphisms. Then $\ch_{\cof(\cx)}$ is a full image
of an accessible functor into $\ck^\to$.
\end{propo}  

\begin{proof}
Given $n<\omega$, let $\cm_n$ be the category whose objects are $(4n+2)$-tuples
$$
(f,g,a_1,\dots,a_n,b_1,\dots,b_n,h_1,\dots,h_n,k_1,\dots,k_n)
$$
of morphisms $f:A\to B$, $g:B\to A$, $a_1,\dots,a_n:A\to A$, $b_1,\dots,b_n:B\to B$, $h_1,\dots,h_n:C(A)\to A$
and $k_1,\dots,k_n:C(B)\to B$. Morphisms are pairs $(u,v)$ of morphisms $u:A\to A'$ and $v:B\to B'$ such that
$f'u=vf$, $g'v=ug$, $uh_i=h'_iC(u)$ and $vk_i=k'_iC(v)$ for $i=1,\dots,n$. Since the cylinder functor is accessible,
$\cm_n$ is an accessible category (use \cite{AR}, 2.67). Let $\overline{\cm}_n$ be the full subcategory of $\cm_n$
such that $h_1\gamma_A=(gf,a_1)$, $h_i\gamma_A=(a_i,a_{i+1})$, $h_n\gamma_n=(a_n,\id_A)$, $k_1\gamma_A=(fg,b_1)$, 
$k_i\gamma_A=(b_i,b_{i+1})$ and $k_n\gamma_n=(b_n,\id_B)$ where $1<i<n$. In the same way as in \cite{AR}, 2.78,
we show that $\overline{\cm}_n$ is accessible as well. We have full embeddings
$$
M_{m,n}:\cm_m\to\cm_n,
$$
for $m<n$, which takes the missing $a_i,b_i,h_i,k_i $ as the identities. By the same reason as above, the union $\cm$ 
of all $\cm_n$'s is an accessible category. Since all $\overline{\cm}_n$'s are accessibly embedded into $\cm$, their
union $\overline{\cm}$ is accessible by \ref{prop2.1}. Let 
$$
F:\overline{\cm}\to\ck^\to
$$
sends each $(4n+2)$-tuple above to $f$. 
This is an accessible functor whose image is $\ch_{\cof(\cx)}$.
\end{proof}

\section{Weak equivalences}
A \textit{model category} is a complete and cocomplete category $\ck$ together with three classes of
morphisms $\cf$, $\cc$ and $\cw$ called \textit{fibrations}, \textit{cofibrations} and \textit{weak 
equivalences} such that
\begin{enumerate}
\item[(1)] $\cw$ has the 2-out-of-3 property, i.e., with any two of $f$, $g$, $gf$ belonging to $\cw$ also
the third morphism belongs to $\cw$, and $\cw$ is closed under retracts in the arrow category $\ck^\to$, and
\item[(2)] $(\cc,\cf\cap\cw)$ and $(\cc\cap\cw,\cf)$ are weak factorization systems.
\end{enumerate}
 
Morphisms from $\cf\cap\cw$ are called \textit{trivial fibrations} while morphisms from $\cc\cap\cw$
\textit{trivial cofibrations}. A \textit{cofibrant replacement functor} $R_c:\ck\to\ck$ is given by 
the (cofibration, trivial fibration)
factorization
$$
0\to R_c(K)\to K
$$
of the unique morphism from $0$ to $K$ while a \textit{fibrant replacement functor} $R_f$ is given by the
(trivial cofibration, fibration) factorization of $K\to 1$. Their composition $R=R_fR_c$ is called 
a \textit{replacement functor}. All three replacements $R_c$, $R_f$ and $R$ can be taken as functors
$\ck^\to\to\ck^\to$ as well. 

A model category $\ck$ is called \textit{combinatorial} provided that $\ck$ is locally presentable and the both 
weak factorization systems $(\cc,\cf\cap\cw)$ and $(\cc\cap\cw,\cf)$ are cofibrantly generated. In a combinatorial
model category, all three replacement functors $R_c$ ,$R_f$ and $R$ are accessible both as functors on $\ck$
and as functors on $\ck^\to$.

\begin{theo}\label {th4.1}
Let $\ck$ be a combinatorial model category. Then $\cw$ is accessible and accessibly embedded into $\ck^\to$.
\end{theo}
\begin{proof}
There is a regular cardinal $\lambda$ such that $\cw$ is closed in $\ck^\to$ under $\lambda$-directed colimits
(see \cite{D}, 7.3). A morphism $f:K\to L$ is a weak equivalence if and only if its replacement $R(f)$ is
a homotopy equivalence. Thus $\cw$ is a pseudopullback
$$
\xymatrix@=4pc{
\ck^\to \ar[r]^{R} & \ck^\to \\
\cw \ar [u]^{} \ar [r]_{} &
\ch_\cc \ar[u]_{}
}
$$
Following \ref{prop3.8} and \ref{lem2.6}, $\cw$ is a full image of an accessible functor $\cm\to\ck^\to$. Thus
it is accessible (see \ref{rem2.3}).
\end{proof}

\begin{theo}\label{th4.2}
Let $\cx$ be a set of morphisms in a locally presentable category $\ck$. Then $\cc=\cof(\cx)$ and $\cw$ make $\ck$ 
a combinatorial model category if and only if
\begin{enumerate}
\item[(1)] $\cw$ has the 2-out-of-3 property and is closed under retracts in $\ck^\to$,
\item[(2)] $(\cx)^\square\subseteq\cw$,
\item[(3)] $\cof(\cx)\cap\cw$ is closed under pushout and transfinite composition, and
\item[(4)] $\cw$ satisfies the solution set-condition at $\cx$.
\end{enumerate}
\end{theo}
\begin{proof}
Sufficiency was shown in \cite{B} and necessity follows from \ref{th4.1} and \cite{AR}, 2.45.
\end{proof}

\begin{coro}\label{cor4.3}
Let $\cx$ be a set of morphisms in a locally presentable category $\ck$ and $\cw_i$, $i\in I$ be a set of classes
of morphisms of $\ck$. Let $\cc=\cof(\cx)$ and $\cw_i$ make $\ck$ a combinatorial model category for each $i\in I$.
Then $\cc$ and $\cap_{i\in I} \cw_i$ make $\ck$ a combinatorial model category.
\end{coro}
\begin{proof}
Since an intersection of a set of accessible and accessibly embedded subcategories of $\ck^\to$ is accessible
(see \cite{AR}, 2.37), the result follows from \ref{th4.1} and \ref{th4.2}.
\end{proof}

\begin{coro}\label{cor4.4}
Let $\cx$ be a set of morphisms in a locally presentable category $\ck$. Assuming Vop\v enka's principle,
$\cc=\cof(\cx)$ and $\cw$ make $\ck$ a combinatorial model category if and only if
\begin{enumerate}
\item[(1)] $\cw$ has the 2-out-of-3 property and is closed under retracts in $\ck^\to$,
\item[(2)] $(\cx)^\square\subseteq\cw$, and
\item[(3)] $\cof(\cx)\cap\cw$ is closed under pushout and transfinite composition.
\end{enumerate}
\end{coro}
\begin{proof}
Vop\v enka's principle implies that each full subcategory of $\ck^\to$ is cone-reflective (see \cite{AR}, 6.7).
Thus condition (4) in \ref{th4.2} is automatic.
\end{proof} 

The statement is equivalent to Vop\v enka's principle as the following example demonstrates.

\begin{exam}\label{ex4.5}
{
\em
Let $\ca$ be a reflective full subcategory of a locally finitely presentable category $\ck$. The corresponding reflector
will be denoted as $R$. Let $\ce$ be the class of all morphisms $f$ such that $R(f)$ is an isomorphisms. Then
$(\ce,\ce^\bot)$ is a factorization system where $\ce^\bot$ consists of morphisms having a unique lifting property
w.r.t. $\ce$ (it means that the diagonal $d$ in \ref{def3.1} is unique). Then $\ce$ satisfies conditions (1) and (3)
from \ref{cor4.4} (see, e.g., \cite{RT1}). 

The weak factorization system $(\ck^\to,\Iso)$ is cofibrantly generated by the set $\cx$ consisting of morphisms $f$ 
having finitely presentable domains and codomains. This can be found in the dual form in \cite{DR}, Proposition 3.1 
(see also \cite{AR2}). In fact, $g\in\cx^\square$ if and only if it is both a pure monomorphism and a pure epimorphism.
Consequently, $\cc=\ck^\to$ and $\cw=\ce$ satisfies all assumptions in \ref{cor4.4}. Thus there is a regular cardinal
$\lambda$ such that $\ce$ is $\lambda$-accessible and closed in $\ck^\to$ under $\lambda$-directed colimits. Let
$\ce_0$ consist of $\lambda$-presentable objects in $\ce$. Since $\ce$ is the closure of $\ce_0$ under $\lambda$-directed
colimits in $\ck^\to$, we have
$$
\ce^\bot=\ce_0^\bot
$$
(see \cite{FR}, 2.2). Since an object $K$ is orthogonal to $\ce$ if and only if its unique morphism $t:K\to 1$
to
the terminal object belongs to $\ce^\bot$, the reflective subcategory $\ca$ is a small-orthogonality class. We use
the fact that $K\in\ca$ if and only if $K\to 1$ belongs to $\ce^\bot$ (see, e.g., \cite{RT}, 3.2). 

We have proved that the statement of \ref{cor4.4} implies that every reflective full subcategory of a locally
finitely presentable category is a small-orthogonality class. But, this implies Vop\v enka's principle (see
\cite{AR}, 6.24.)
}
\end{exam}

\begin{coro}\label{cor4.6}
Let $\cx$ be a set of morphisms in a locally presentable category $\ck$ and $\cw_i$, $i\in I$ be a collection of classes
of morphisms of $\ck$. Let $\cc=\cof(\cx)$ and $\cw_i$ make $\ck$ a combinatorial model category for each $i\in I$.
Then $\cc$ and $\cap_{i\in I} \cw_i$ make $\ck$ a combinatorial model category.
\end{coro}
\begin{proof}
It follows from \ref{cor4.4}.
\end{proof}

\begin{rem}\label{rem4.7}
{
\em
Let $\cx$ be a set of morphisms of a locally presentable category. Since $\cc=\cx^\to$ and $\ck^\to$ always form
a combinatorial model category, $\cw=\ck^\to$ is the largest $\cw$ with this property. Assuming Vop\v enka's principle, 
these classes $\cw$ form a (large) complete lattice. Thus there is the smallest class $\cw_\cx$ with this property.
Moreover, it can be constructed as the closure of $\cx^\square$ under properties (1) and (3) from \ref{th4.2}.
This construction was introduced in \cite{RT} where we called $\cw_\cx$ \textit{left determined}. Independently,
this construction was considered by D.-C. Cisinski \cite{C} who proved, without any set theory, that one 
gets a combinatorial model category in the special case when $\ck$ is a Grothendieck topos and $\cof(I)$
is the class of all monomorphisms (the latter class is always cofibrantly generated in this case). 

In 2002, J. H. Smith informed me that he is able to prove that the $\cw$'s above form a small complete lattice
without any set theory. Thus it always has the smallest element. Under Vop\v enka's principle, this is our
left-determined $\cw_\cx$. I do not know what happens in general and, in particular, I have not been able to
prove the Smith's claim.
}
\end{rem}

\end{document}